\documentclass[]{svmult} 

\usepackage[utf8]{inputenc} 

\RequirePackage{amsfonts,amssymb,stmaryrd}
\RequirePackage{mathtools, xcolor}
\usepackage{breqn, siunitx, url}
\RequirePackage{bm, graphicx, xspace}
\usepackage{commands}

\usepackage{array} % for better arrays (eg matrices) in maths
\usepackage{subfig} % make it possible to include more than one captioned figure/table in a single float

\begin{document}

\title*{Finite volume discretisation of fracture deformation in thermo-poroelastic media}

\author{Ivar Stefansson \and Inga Berre \and Eirik Keilegavlen}
\institute{Ivar Stefansson \email{ivar.stefansson@uib.no} \and Inga Berre  \and Eirik Keilegavlen  \at University of Bergen, Bergen, Norway}
\date{} % Activate to display a given date or no date (if empty),
         % otherwise the current date is printed 
\maketitle

\abstract*{This paper presents a model where thermo-hydro-mechanical processes are coupled to a deformation model for preexisting fractures. 
The model is formulated within a discrete-fracture-matrix framework where the rock matrix and the fractures are considered as individual subdomains, and interaction between them takes place on the matrix-fracture interfaces.
A finite volume discretisation implemented in the simulation toolbox PorePy is presented and applied in a simulation showcasing the effects of the different mechanisms on fracture deformation governed by contact mechanics, as well as their different timescales. 
}

\abstract{This paper presents a model where thermo-hydro-mechanical processes are coupled to a deformation model for preexisting fractures. 
The model is formulated within a discrete-fracture-matrix framework where the rock matrix and the fractures are considered as individual subdomains, and interaction between them takes place on the matrix-fracture interfaces.
A finite volume discretisation implemented in the simulation toolbox PorePy is presented and applied in a simulation showcasing the effects of the different mechanisms on fracture deformation governed by contact mechanics, as well as their different timescales. 
}

\section{Introduction}\label{sec:introduction}
We consider the simulation of fully coupled thermo-hydro-mechanical (THM)
dynamics in fractured porous media, where the fractures can undergo sliding
if the shear forces on the fracture planes are sufficient to overcome 
frictional resistance.
These processes are highly relevant for several subsurface applications,
including geothermal energy extraction, storage of CO$_2$ and energy and
groundwater management. 
Our simulation approach is based on three main ingredients:
First, conservation of mass, energy and momentum is preserved under
discretisation by the employment of a fully coupled finite volume approach for
the governing equations. 
Second, the network of fractures, which act as main conduits for fluid
flow and energy transport, is explicitly represented in the simulation model.
Specifically, the fractures are represented as lower-dimensional manifolds
that are embedded in the host porous medium, thus the simulation model is
defined on a mixed-dimensional geometry.
Third, the sliding of fractures is modelled as a frictional contact problem,
which is solved by an active set approach.
The discretisation of the contact problem benefits from the finite
volume approach, which directly provides discrete representations of displacements
as well as of mechanical, fluid and thermal forces  on the fracture surfaces.
Furthermore, the explicit degrees of freedom for fluid pressure inside the fractures
allow us to capture the critical interplay between elevated fluid pressures and
fracture deformation.

\section{Model}\label{sec:model}
%\subsection{Mixed-dimensional geometry} 
%The host porous medium and the fractures are together represented as a
We consider a mixed-dimensional geometry which is decomposed in
subdomains of different dimensions representing the host porous medium and
the lower-dimensional planar fractures, and separated by interfaces, see \cite{keilegavlen2019porepy} for details.
Variables and governing equations are defined on subdomains and interfaces,
with full flexibility to vary the type and number of variables and equations
between geometric objects.
The framework accommodates heterogeneous and multiphysics models, with a natural
treatment of modelling and discretisation of mixed-dimensional problems.

\begin{figure}[bp]
\centering
\includegraphics[width=.45\textwidth]{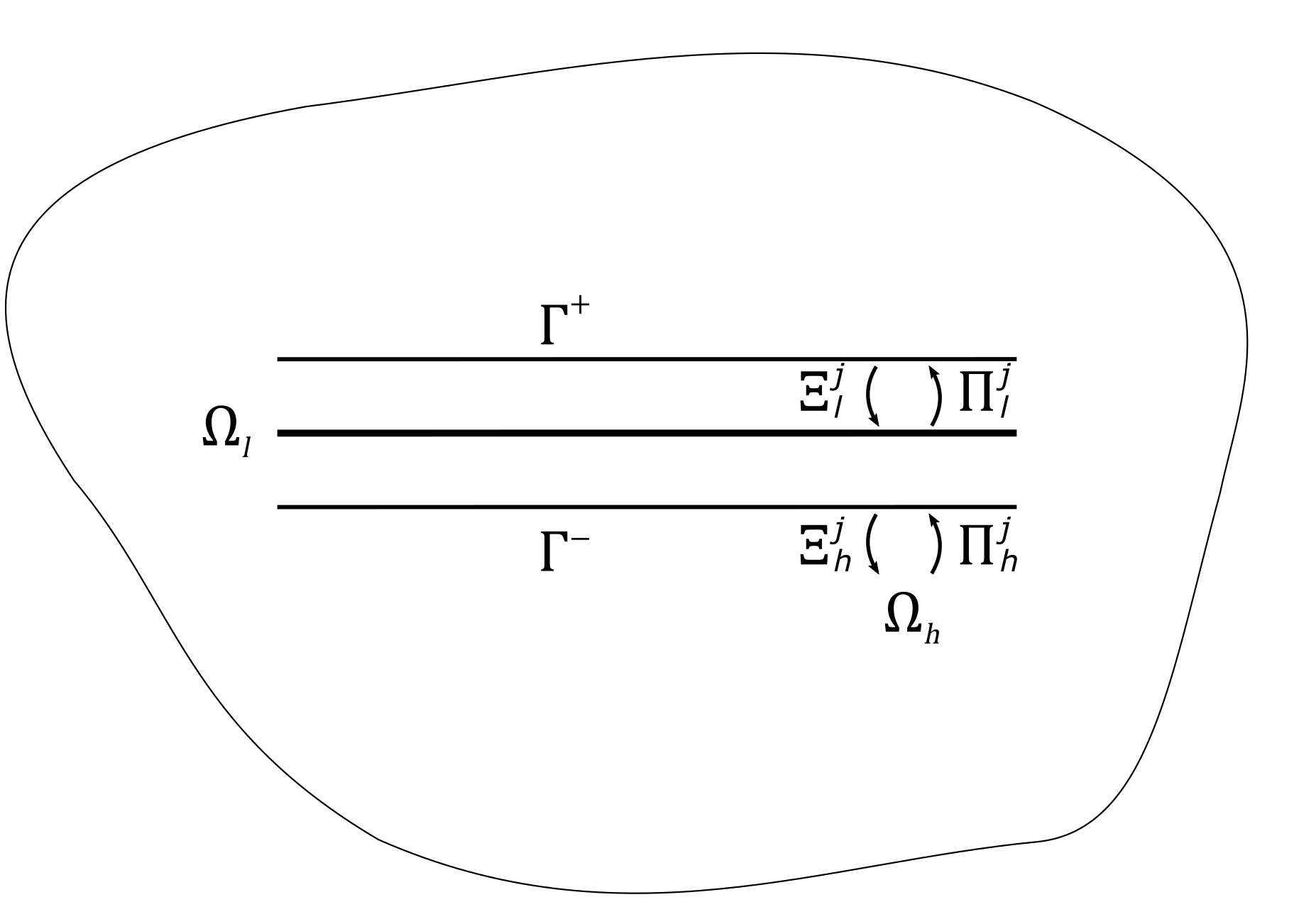}
\includegraphics[width=.54\textwidth]{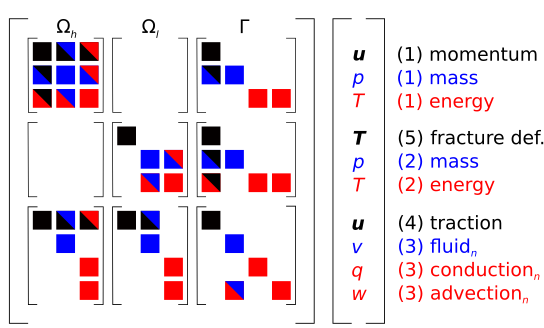}
\caption{Schematic representation of a fracture \domain[l] and a
	matrix subdomain \domain[h] to the left. The two subdomains are
	separated by the interface, whose two sides are denoted by
	$\interface^+$ and $\interface^-$. Also shown are the projection
	operators used for transfer of variables between the subdomains and the
	interface.
	To the right, we show the block structure of the matrix resulting from a 
	fully implicit discretisation of Eqs. \ref{eq:conservation_equations_matrix} 
	through \ref{eq:fracture_conditions}.
}
\label{fig:fracture_schematic_projections_and_interface_pm}
\end{figure}

We denote a subdomain by \domain[i] and
its boundary by \boundary[i], and identify the variables defined within it by
subscript $i$.  Where convenient, we will denote the higher-dimensional matrix
domain by \domain[h] and lower-dimensional fracture domains by \domain[l],
as indicated in Fig. \ref{fig:fracture_schematic_projections_and_interface_pm}.

\boundary[i] may be divided into the external boundary $\boundary[i]^{e}$ and
the internal fracture boundary $\boundary[i]^{f}$, which coincides geometrically
with both the immersed fracture domain \domain [l] and the interface between
\domain[h] and \domain[l].  This interface is denoted by \interface[j], and the
associated variables identified by subscript $j$.  The two
sides of the fracture are denoted by $+$ and $-$, as shown in Fig.
\ref{fig:fracture_schematic_projections_and_interface_pm}.  
Projection of variables from the interface to the subdomains is
performed by \projectFromInterface{h}{j} and $\projectFromInterface{l}{j}$,
respectively, whereas \projectToInterface{h}{j} and \projectToInterface{l}{j}
project from the subdomains to the interface.

%\subsection{Governing equations and constitutive laws}
For \domain[h], the primary variables are displacement \displacement, 
fluid pressure \pressure and temperature \temperature. 
The fluid flux and the advective and conductive heat fluxes are denoted
by \fluidFlux, \advectiveHeatFlux and \conductiveHeatFlux, respectively.
Assuming all external source and sink terms to be zero, 
conservation of momentum, mass and energy in \domain[h] can be defined as \cite{salimzadeh2018thm}

\begin{gather}\label{eq:conservation_equations_matrix}
\begin{aligned}
\divergence  [ \frac{\stiffness}{2}(\nabla\displacement + \nabla\displacement^T) - \biotAlpha \pressure \identity - \thermalExpansion[s] \bulkModulus (\temperature - \temperature[0]) \identity]
&= 0, \\
  \left(\porosity \compressibility+\frac{\biotAlpha-\porosity}{\bulkModulus}\right) \dd{\pressure}{\timet} +  \biotAlpha \dd{(\divergence\displacement)}{\timet} - \porosity \thermalExpansion[f] \dd{\temperature}{\timet} - \divergence \frac{\permeability}{\viscosity} \nabla \pressure &= 0,\\
\density \heatCapacity \dd{\temperature}{\timet} +  \thermalExpansion[s] \bulkModulus \temperature[0] \dd{(\divergence\displacement)}{\timet} - \porosity \thermalExpansion[f] \temperature[0] \dd{\pressure}{\timet} + \density[f] \heatCapacity[f] \fluidFlux \cdot \nabla \temperature - \divergence \heatConductivity \nabla \temperature&=0.
\end{aligned}
\end{gather}
Similarly, conservation of mass and energy in 
\domain[l] and fracture intersection domains \domain[x] is given by
\begin{gather}\label{eq:conservation_equations_fractures}
\begin{aligned}
 \specificVolume \compressibility \dd{\pressure}{\timet} +  \dd{\specificVolume}{\timet}- \specificVolume \thermalExpansion[f] \dd{\temperature}{\timet} - \divergence \frac{\permeability}{\viscosity} \nabla \pressure  &=\projectFromInterface{i}{j}\interfaceFluidFlux, \\
\specificVolume\density \heatCapacity \dd{\temperature}{\timet} - \specificVolume \thermalExpansion[f] \temperature[0] \dd{\pressure}{\timet} + \density[f] \heatCapacity[f] \fluidFlux \cdot \nabla \temperature - \divergence \heatConductivity \nabla \temperature&= \projectFromInterface{i}{j}\interfaceConductiveHeatFlux + \projectFromInterface{i}{j}\interfaceAdvectiveHeatFlux.
\end{aligned}
\end{gather}
Here, the specific volume \specificVolume accounts
for the extension in the collapsed dimension(s), while subscript 0 denotes 
the initial value and $f$ and $s$ indicate fluid and solid parameters, respectively.

Denoting the trace operator by \trace{\cdot}, the conditions on \interface[j]
 are the three flux relationships and the traction balance:
\begin{gather}\label{eq:interface_conditions}
\begin{aligned}
      \interfaceFluidFlux &= - \permeability[j] (\projectToInterface{l}{j} \pressure[l] - \projectToInterface{h}{j} \trace{\pressure[h]}), \\
            \interfaceConductiveHeatFlux &= - \heatConductivity[j] (\projectToInterface{l}{j} \temperature[l] - \projectToInterface{h}{j} \trace{\temperature[h]}), \\
            \interfaceAdvectiveHeatFlux &= \left\{ \begin{array}{ l l }
    \density[f] \heatCapacity[f] \interfaceFluidFlux \projectToInterface{h}{j}\trace{\temperature[h]} & \quad\qquad \text{ if } \interfaceFluidFlux>0,  \\
    \density[f] \heatCapacity[f] \interfaceFluidFlux \projectToInterface{l}{j}\temperature[l] & \qquad\quad \text{ if } \interfaceFluidFlux \leq 0. 
  \end{array} \right. \\            
    \end{aligned}
\end{gather}
 Eqs. \ref{eq:conservation_equations_matrix} through \ref{eq:interface_conditions} are complemented by the internal boundary conditions
$	\trace{\displacement[h]} = \projectFromInterface{h}{j}\displacement[j] $,
$	\fluidFlux[h]\cdot\normal = \projectFromInterface{h}{j}\interfaceFluidFlux[j] $,
$	\conductiveHeatFlux[h]\cdot \normal = \projectFromInterface{h}{j}\interfaceConductiveHeatFlux[j]$ and
$	\advectiveHeatFlux[h]\cdot \normal = \projectFromInterface{h}{j}\interfaceAdvectiveHeatFlux[j] \text{ on }\boundary[h]^{f}$,
 and standard Dirichlet and Neumann conditions on the external boundaries.

The fracture deformation is described by relations between the contact traction on the
fracture surface, \traction, and the jump in displacement over the fracture.
Denoting  the displacements on the two sides of the interface by 
$\displacement[j]^+$ and $\displacement[j]^-$, the displacement jump is
$\jump{	\displacement[j]} =\projectFromInterface{l}{j}\left(\displacement[j]^+-\displacement[j]^-\right)$,
and \jump{\increment\displacement} denotes its increment.
$\jump{\displacement}$ is also related to the aperture \aperture and specific volume: for \domain[l], 
we set $\specificVolume=\aperture=\jump{\displacement}_n$. For \domain [x] 
we use the mean and product of \aperture for the adjacent cells of the intersecting fractures when computing 
\aperture and \specificVolume, respectively.

Since the fracture deformation depends on the traction caused by the \textit{contact} 
between the two surfaces, we subtract the contribution from the pressure \pressure[l] on the fracture surfaces.
Thus, the interface tractions on the two fracture surfaces and the traction balance are
\begin{gather}\label{eq:fracture_interface_traction}
  \begin{aligned}
    \traction[j]^+ &=  \projectToInterface{h}{j} \stress \cdot \normal|_{\boundary[h]^+}, \\
    \traction[j]^- &=\projectToInterface{h}{j} \stress \cdot \normal|_{\boundary[h]^-}, \\
    \projectToInterface{l}{j}(\traction+\pressure[l])&= \traction[j]^+ = - \traction[j]^-.
  \end{aligned}
\end{gather}

Denoting tangential and normal components of vectors on the fracture by subscripts \tangential
 and $n$, respectively, the fracture deformation for \domain[l] is governed by three non-penetration relations and three friction law constraints:

\begin{gather}\label{eq:fracture_conditions}
\begin{aligned}
\begin{array}{ r l r l}
 \jump{\displacement}_n   &\leq 0  &  
     \norm{\traction[\tau]}&\leq -F\normalTraction\\
 \normalTraction \jump{\displacement}_n &= 0    &      
     \norm{\traction[\tau]} &< -\frictionCoefficient \normalTraction \,\,\rightarrow \,\,\jump{\increment\displacement}_{\tau} = 0 \\
 \normalTraction &\leq 0 &  
      \qquad\norm{\traction[\tau]}&= -\frictionCoefficient \normalTraction \,\,\rightarrow \,\,\exists \,\zeta \in \mathbb{R^-}: \traction[\tau] = \zeta \jump{\increment\displacement}_{\tau}.
\end{array}
\end{aligned}
\end{gather}
Further detail on the fracture deformation is found in \cite{berge2019hmdfm, hueeber2008}, whereas 
the definition of all parameters may be found in the repository at \cite{stefansson2019runscripts}.

\section{Discretisation}\label{sec:discretisation}
Applying Implicit Euler for the temporal discretisation, the scalar conservation Eqs.
\ref{eq:conservation_equations_matrix} and \ref{eq:conservation_equations_fractures}
are discretised using a Multi-Point Flux Approximation \cite{aavatsmark2002mpfa} for the conductive terms and a first order upwind scheme for the advective term.  
The momentuum conservation equation and the \divergence \displacement terms in the scalar conservation laws are discretised using the FV scheme introduced in \cite{nordbotten2016biot}. 
The scheme, termed Multi-Point Stress Approximation (MPSA), is based on local momentum conservation and is formulated in terms of discrete cell centred pressures and displacement unknowns.
Originally developed for the pure hydro-mechanical problem, the coupled discretisation approach can readily be extended to the THM case.

Thanks to the structure provided by the mixed-dimensional framework, the discretisation of the coupling fluxes of Eq. \ref{eq:interface_conditions} consists of two simple tasks. 
Discrete projection operators transfer variables from higher-dimensional faces and lower-dimensional cells to the interface cells.
The interface fluxes are then discretised directly using the projected variables, see \cite{keilegavlen2019porepy}.

The traction balance and fracture deformation relations of Eqs. 
\ref{eq:fracture_interface_traction} and \ref{eq:fracture_conditions} are
formulated in terms of displacement and traction on the fracture surfaces.  The
former is included as a primary interface variable, and thus directly
available.  
While the latter is not a primary variable, the FV framework
is formulated in terms of face tractions, implying that in the discrete setting, 
the surface traction can be reconstructed from the primary
variables.
Specifically, we apply available discretisation operators to get contributions to the stress
from displacements, pressures and temperatures in \domain[h], the interface variable \displacement[j], 
and conditions on external boundaries.

For the discretisation of the fracture deformation relations we first
reformulate Eq. \ref{eq:fracture_conditions} as two nonlinear complementary
functions and compute their derivatives. A semismooth Newton scheme is applied
on the basis of the three sets
%EK: Respectively
\begin{gather}\label{eq:deformation_state_sets}
    \begin{aligned}
    \inactiveSet_n &= \left\{\frictionBound \leq 0 \right\} \\
    \inactiveSet_{\tangential}&= \left\{||-\traction[\tangential]+\cNum\jump{\increment\displacement}_{\tangential}||<\frictionBound \right\} \\
    \activeSet &=\left\{||-\traction[\tangential]+\cNum\jump{\increment\displacement}_{\tangential}||\geq \frictionBound >0  \right\},
    \end{aligned}
  \end{gather}
which correspond to fracture cells which are open, sticking and sliding, respectively.
 \cNum denotes a numerical parameter and $\frictionBound=\frictionCoefficient\left(-\normalTraction+\jump{\displacement}_{n} \right)$ 
 is the friction bound. 
For each fracture cell \cell, this results in the following cell-wise constraints when computing iterate $k+1$ from the current iterate $k$:
\begin{gather}\label{eq:deformation_constraints}
    \begin{aligned}
	   \begin{array} {r l l}	   
	   \traction^{k+1} &= \vectorFont{0} & \cell \in \inactiveSet_n \\	
	   \jump{\displacement^{k+1}}_n &=0 & \cell \in \inactiveSet_{\tangential} \cup \activeSet \\
	   \jump{\increment\displacement^{k+1}}_{\tangential} + \frictionCoefficient\jump{\increment\displacement^{k}}_{\tangential}/\frictionBound^{k} \normalTraction^{k+1} &= \jump{\increment\displacement^{k}}_{\tangential}& \cell \in \inactiveSet_{\tangential} \\   
	   \traction[\tangential]^{k+1} + L^{k}\jump{\increment\displacement^{k+1}}_{\tangential} + \frictionCoefficient \fractureV^{k} \normalTraction^{k+1} &= \fractureR^{k} + \frictionBound^{k}\fractureV^{k}\qquad\qquad & \cell \in \activeSet.
	   \end{array}
    \end{aligned}
\end{gather}
The coefficients \fractureL, \fractureV and \fractureR are functions of
$\jump{\increment\displacement^{k}}_{\tangential}$, $\jump{\displacement^{k}}_{n}$
and $\traction^{k}$, and can thus be computed from the previous iterate
and time step. For further details of the discretisation and implementation of
the fracture deformation equations, we refer to \cite{berge2019hmdfm}. 

In terms of implementation, we mention that the mixed-dimensional framework
allows us to discretise each term for each subdomain or interface independently.
We may thereby break the highly complex task of discretising the contact conditions 
with a coupled THM stress down in manageable tasks, and also reuse discretisations
for different subdomains. 
For the global discretisation matrix, this manifests as a two-level block structure.  The first level
has the subdomains on the diagonal and interfaces on the off-diagonals. The second level corresponds
to the primary variables, and has coupling terms between different variables on
the off-diagonals.

The model is implemented for two- and three-dimensional problems in the
open source simulation toolbox PorePy presented 
in \cite{keilegavlen2019porepy}, and run scripts for the example simulation 
presented in the following section may be found in the repository \cite{stefansson2019runscripts}.
The simplicial spatial grid is constructed such that the lower-dimensional cells coincide with 
higher-dimensional faces through a back-end to Gmsh \cite{gmsh2009}.

\section{Results}\label{sec:results}

\begin{figure}[tp]
\centering
\includegraphics[width=.8\textwidth]{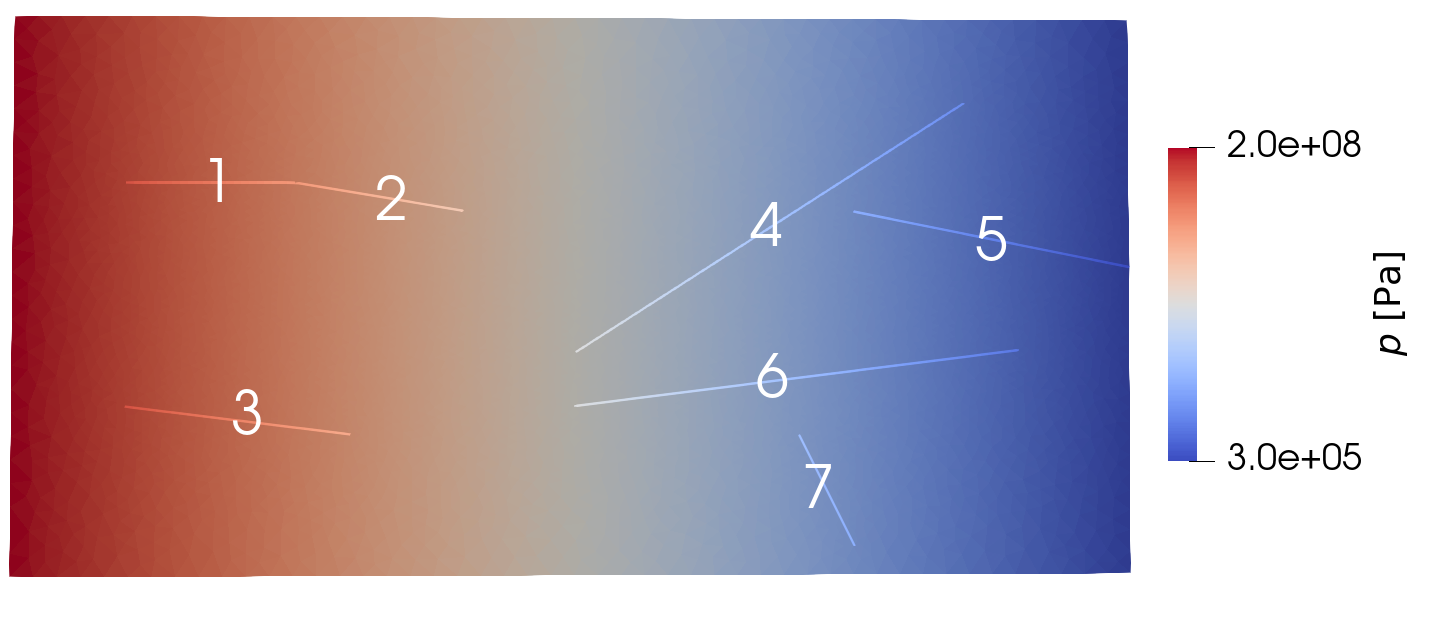}
\caption{The equilibrated pressure state at the end of phase II.
}
\label{fig:pressure_field}
\end{figure}

To demonstrate the applicability of the model and discretisation, we present 
simulations of THM and fracture deformation effects for a 2d domain containing 
seven fractures, see Fig. \ref{fig:pressure_field} for the geometry and
numbering of the fractures. The setup is designed to
expose the method to a wide range of physical driving forces, and thus probe the
stability and performance of the simulation model.
% EK: Strukturer
% EK: Siste setning i intro skal oppsummere denne seksjonen

Starting out from a homogeneous initial state for all primary variables, the
simulation consists of four phases, where we study the effect of sequentially
adding different driving forces.  In phase I, the deformation is caused
by a boundary displacement of
$(0.002,-0.005)^T$ \si{\meter} applied at the top.  To
allow the system to reach equilibrium, this phase lasts from
$\timet=\SI{-10000}{\second}$ to $\timet=\SI{0}{\second}$.  In phase II, a
pressure gradient is applied from left to right. At the end of the phase, at
\timet=\SI{0.02}{\second}, the pressure has virtually reached a steady state,
see Fig. \ref{fig:pressure_field}.  In phase III, we reduce the temperature 
at the left boundary of from \SI{0}{\celsius} to \SI{-100}{\celsius} and
in phase IV we increase it to \SI{100}{\celsius}, thus exploring both thermal
expansion and contraction.  At the end of each of the two last phases, at
$\timet=\SI{2.5}{\second}$ and $\timet=\SI{5}{\second}$, the domain has has reached
a close to uniform temperature.

\begin{figure}[tp]
\centering
\includegraphics[width=.48\textwidth]{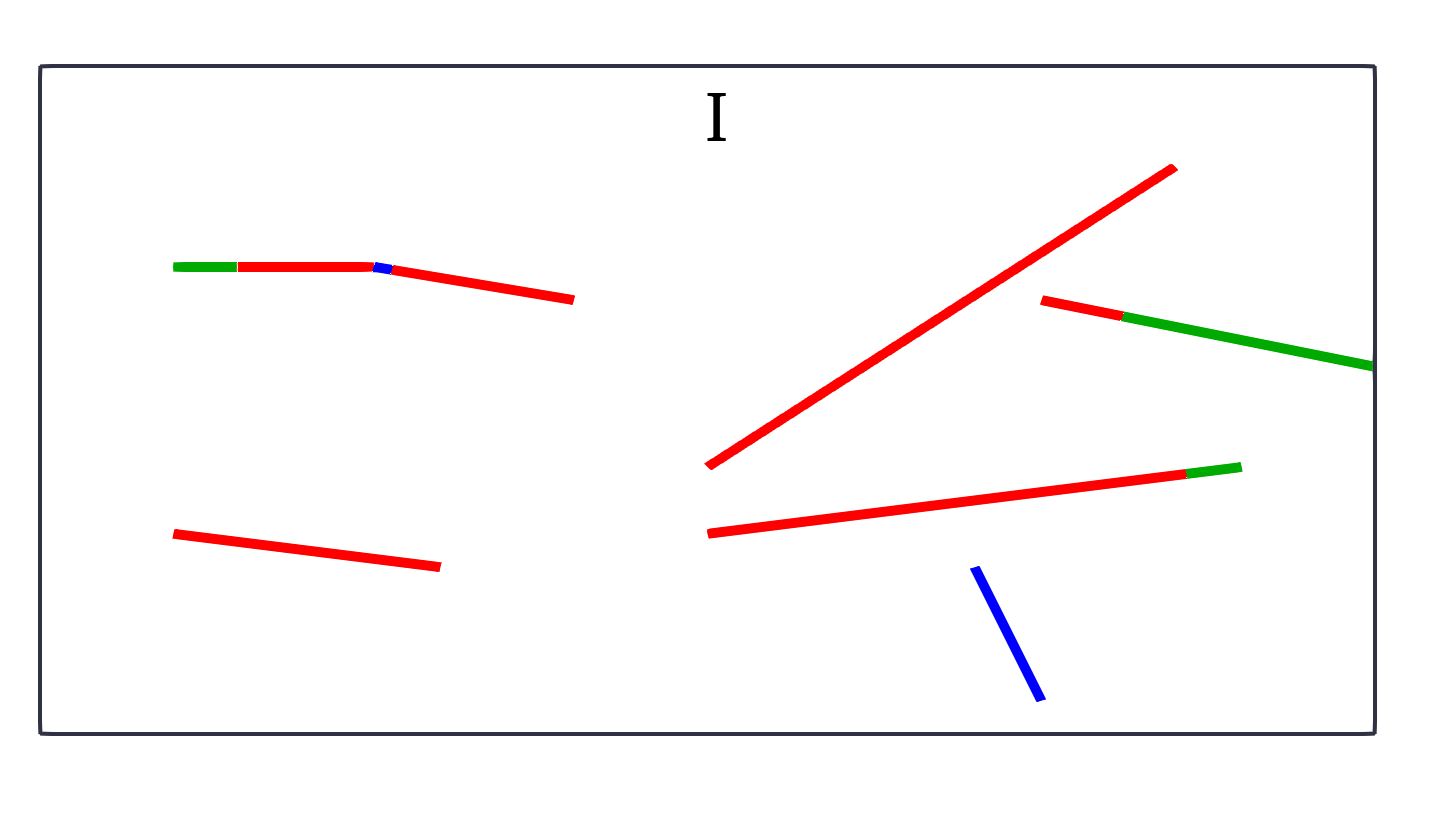}
\includegraphics[width=.48\textwidth]{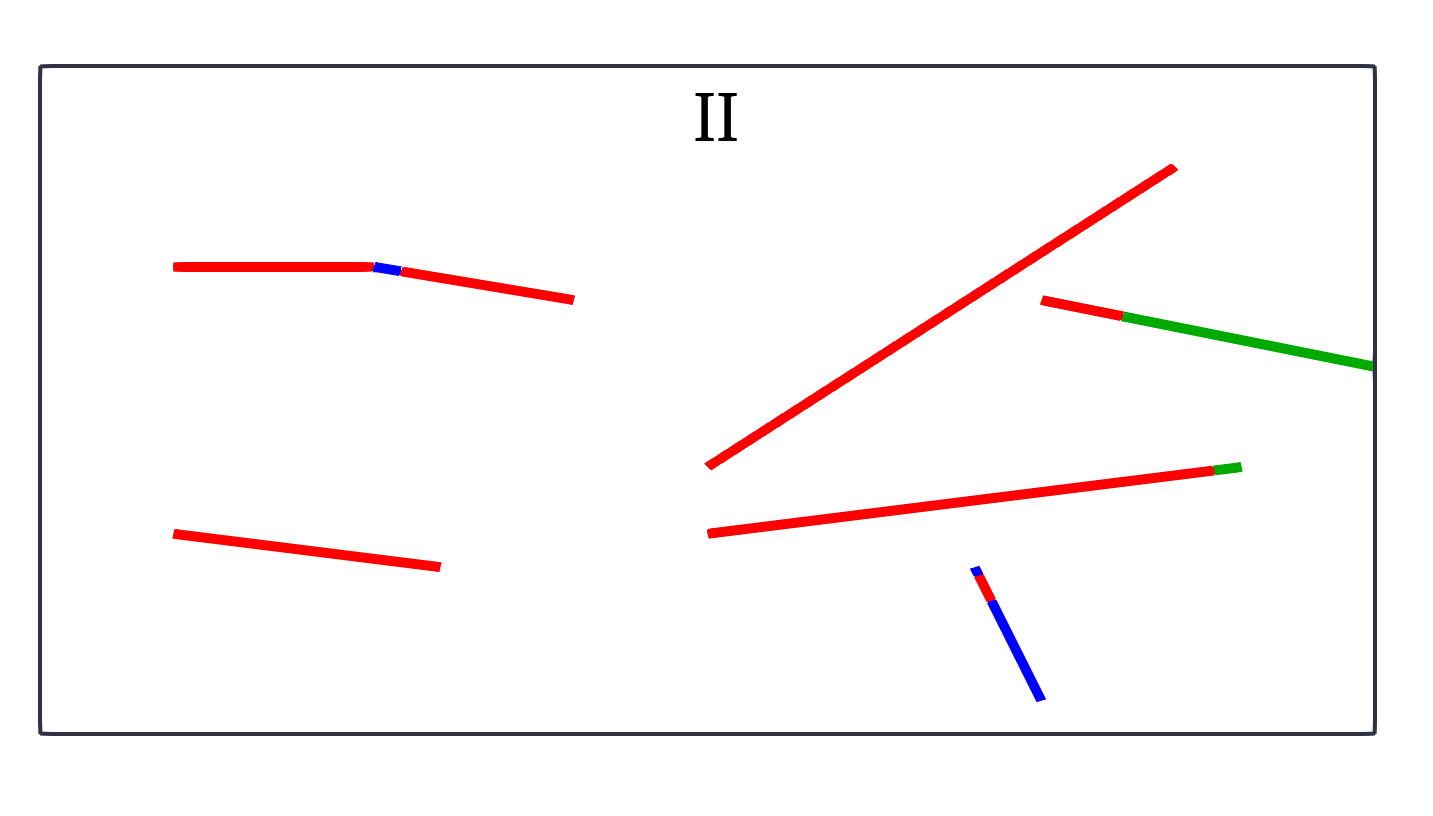}
\includegraphics[width=.48\textwidth]{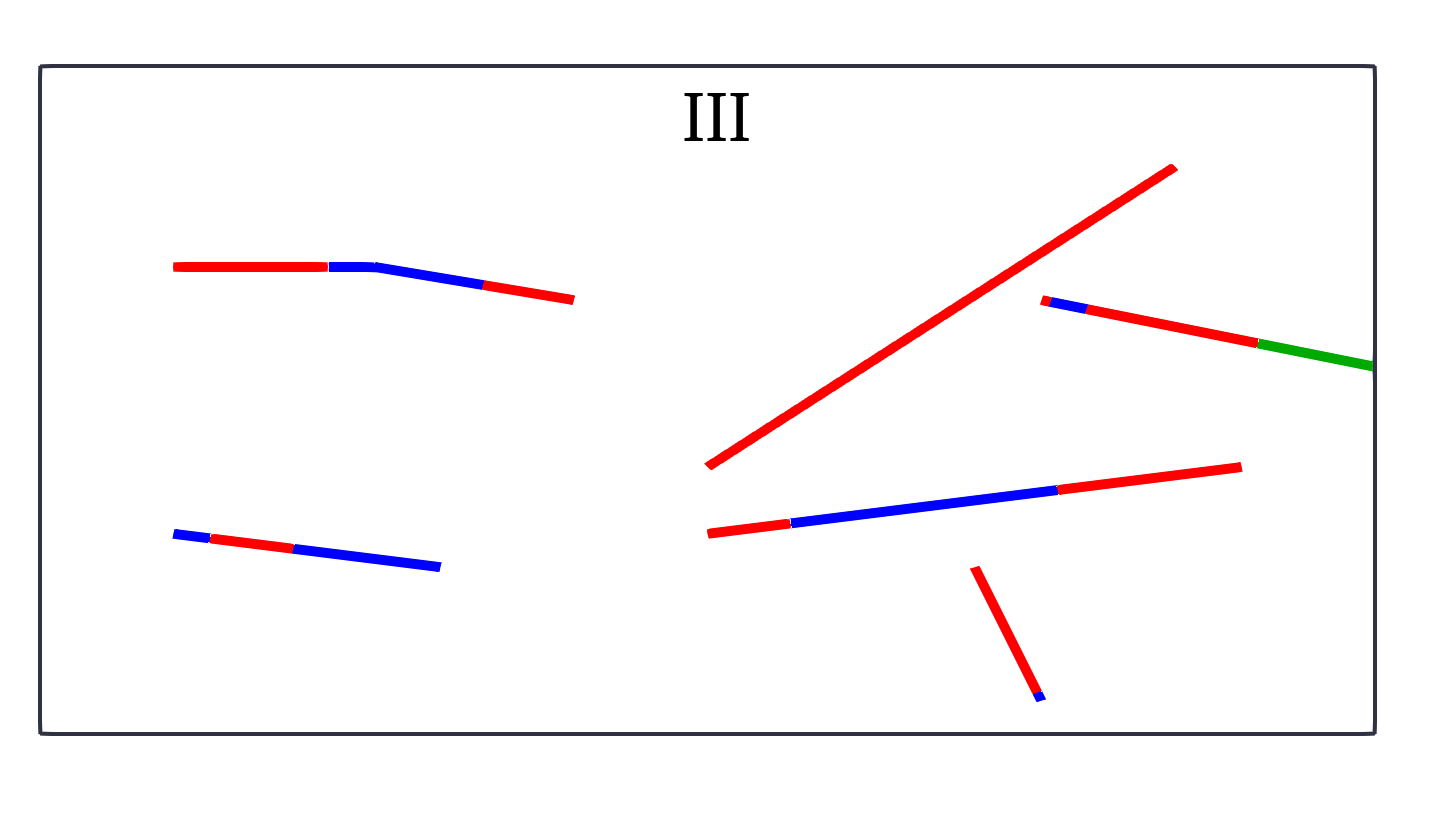}
\includegraphics[width=.48\textwidth]{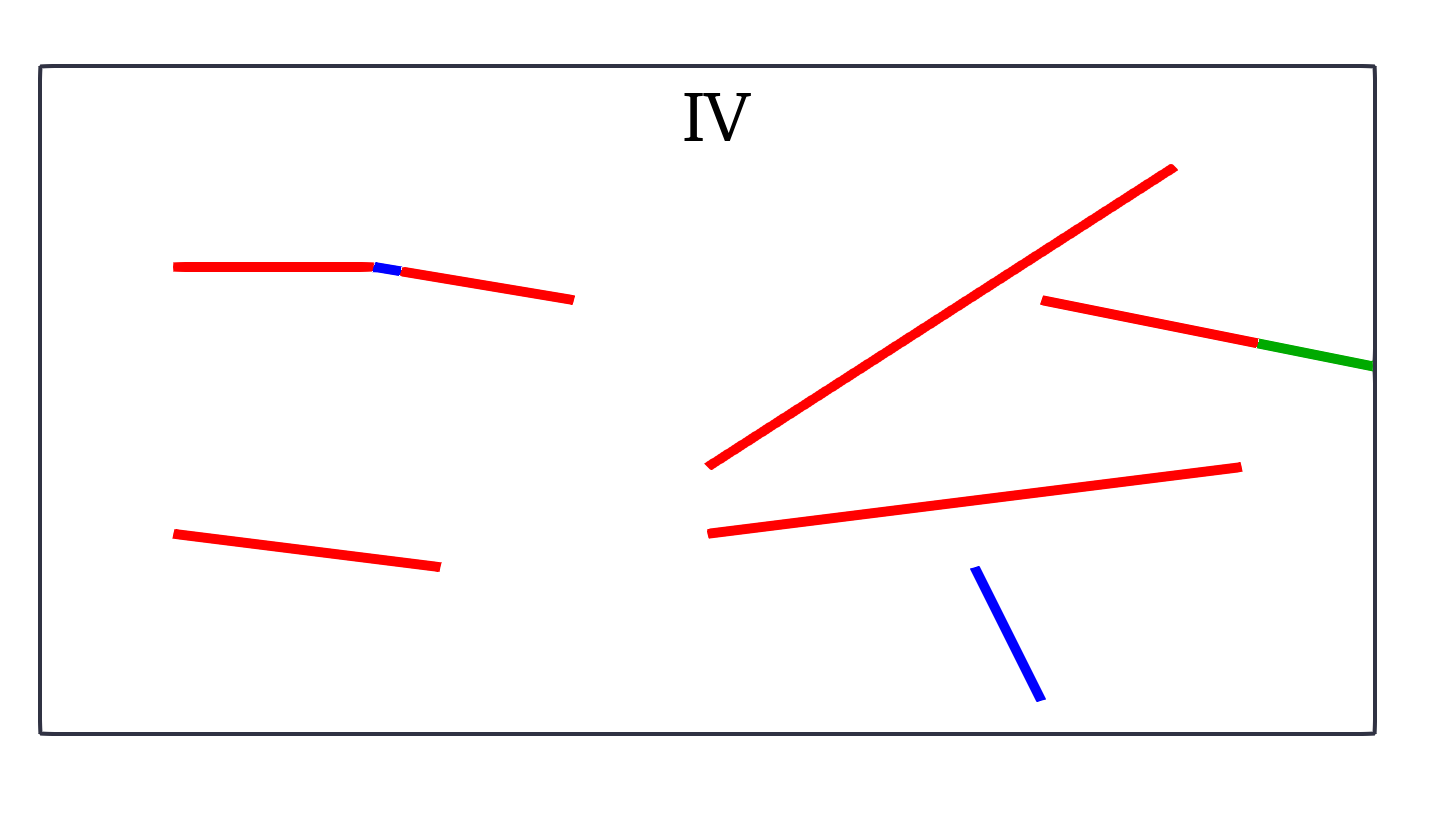} 
\includegraphics[width=.98\textwidth]{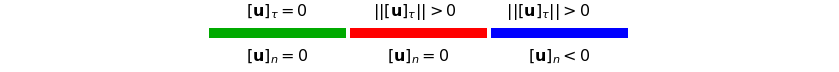} 
\caption{The deformation state of the fractures at the end of the four phases.}
\label{fig:fracture_states}
\end{figure}

\begin{figure}[tp]
\centering
\includegraphics[width=.98\textwidth]{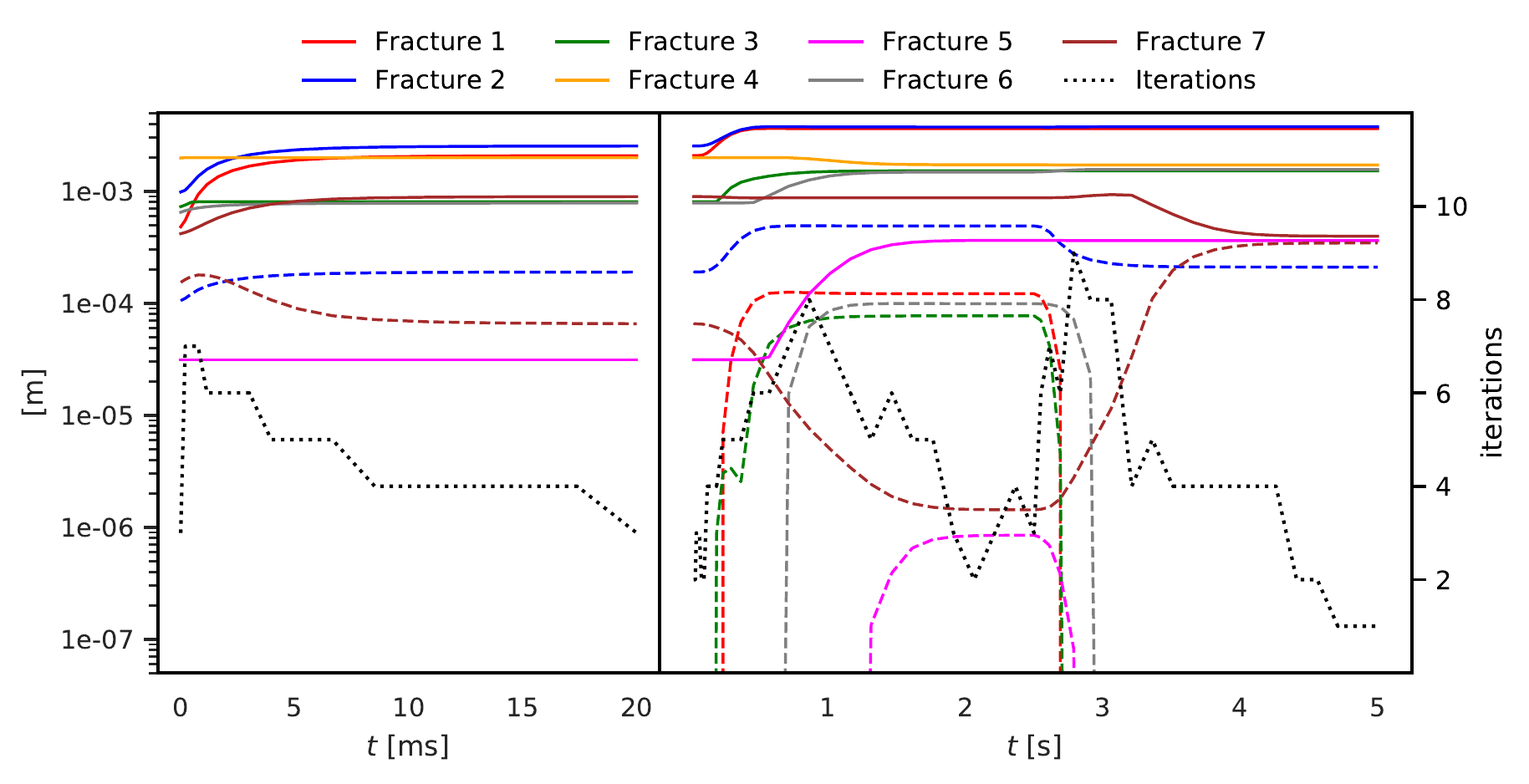} 
\caption{Norm of displacement jumps at each fracture and number of Newton iterations for each time step.
Phase II is shown to the left, with the state at the end of phase I shown at $\timet=0$,
 and phases III and IV to the right.
 Normal jumps are shown in dashed lines and tangential jumps in solid lines. 
}
\label{fig:u_jump_vs_time}
\end{figure}

For the end of each of the simulation phases, Fig. \ref{fig:fracture_states}
shows the deformation state, i.e. whether $\jump{\displacement}_n$ and 
$\jump{\displacement}_\tangential$ are nonzero for each fracture cell.  These
show that the state changes for all fractures, and that for each phase,
at least three fractures change their state.  Fig. \ref{fig:u_jump_vs_time}
shows the norm of $\jump{\displacement}_n$ and
$\jump{\displacement}_\tangential$ on each fracture throughout the simulation. 

Because of the time
scale difference, the pressure phase is shown in the left plot and the
temperature phase in the right one.  The former shows gradual and moderate
deformation for the fractures which have nonzero jumps at the end of phase I,
and onset of sliding for fracture 3.  The latter shows more complex
deformation, displaying non-monotone jump evolution for several fractures, e.g.
fracture 3 first opening and subsequently closing, and fracture 7 undergoing the
reverse process.

Fig. \ref{fig:u_jump_vs_time} also displays the number of Newton iterations
required for convergence for each of the time steps.  While fairly stable
results are demonstrated, some increase may be observed whenever the fracture
state changes more markedly, e.g. when the cooling is introduced at the onset of phase IV.  

\section{Conclusion}\label{sec:conclusion}
A FV framework for simulation of thermo-poroelasticity
fully coupled to fracture deformation is presented.
The contact mechanics problem is naturally discretised
in terms of displacement and traction at the fracture faces, exploiting
the availability of these in the FV formulation of the THM problem.
A numerical example exhibiting complex THM interactions and fracture
dynamics demonstrates both the range of the processes captured by the
model, and its applicability for challenging problems.

\end{document}